\documentclass[sn-basic]{sn-jnl}

\usepackage{ulem}
\usepackage{graphicx}%
\usepackage{multirow}%
\usepackage{amsmath,amssymb,amsfonts}%
\usepackage{amsthm}%
\usepackage{mathrsfs}%
\usepackage[title]{appendix}%
\usepackage{xcolor}%
\usepackage{textcomp}%
\usepackage{manyfoot}%
\usepackage{booktabs}%
\usepackage{listings}%
\usepackage{comment}
\usepackage{algorithm}
\usepackage[noend]{algorithmic}

\theoremstyle{plain}%
\newtheorem{theorem}{Theorem}
\theoremstyle{remark}%

\theoremstyle{definition}%
\newtheorem{definition}{Definition}%

\begin{document}

\title[Stability of PT]{Stability of the persistence transformation}


\author*[1]{\fnm{Gideon} \sur{Klaila}}\email{klailag@uni-bremen.de}
\equalcont{These authors contributed equally to this work.}

\author[1]{\fnm{Anastasios} \sur{Stefanou}}\email{stefanouanastasios@gmail.com}
\equalcont{These authors contributed equally to this work.}

\author[1]{\fnm{Lena} \sur{Ranke}}\email{lranke@uni-bremen.de}
\equalcont{These authors contributed equally to this work.}


\affil[1]{\orgdiv{Insitute for Algebra, Geometry, Topology and its Applications (ALTA), Department of Mathematics}, \orgname{University of Bremen}, \orgaddress{\street{Bibliothekstraße 1}, \postcode{28359}, \city{Bremen}, \country{Germany}}}


\abstract{In this paper, we introduce the persistence transformation, a novel methodology in Topological Data Analysis (TDA) for applications in time series data which can be obtained in various areas such as science, politics, economy, healthcare, engineering, and beyond. This approach captures the enduring presence or `persistence' of signal peaks in time series data arising from Morse functions while preserving their positional information. Through rigorous analysis, we demonstrate that the proposed persistence transformation exhibits stability and outperforms the persistent diagram of Morse functions (with respect to filtration, e.g., the upper levelset filtration). Moreover, we present a modified version of the persistence transformation, termed the reduced persistence transformation, which retains stability while enjoying dimensionality reduction in the data. Consequently, the reduced persistence transformation yields faster computational results for subsequent tasks, such as classification, albeit at the cost of reduced overall accuracy compared to the persistence transformation. However, the reduced persistence transformation finds relevance in specific domains, e.g., MALDI-Imaging, where positional information is of greater significance than the overall signal height. Finally, we provide a conceptual outline for extending the persistence diagram to accommodate higher-dimensional input while assessing its stability under these modifications.}

\keywords{Topological Data Analysis, Time Series, Persistence Homology, Persistence Diagram, Morse functions, Filtrations}



\maketitle

\section{Introduction}
\subsection{Motivation and related work}
Topological Data Analysis (TDA) stands at the intersection of mathematics, data science, and various application domains, promising a unique perspective on understanding complex datasets (e.g., see \cite{carlsson2009topology}, \cite{Edelsbrunner2010}). TDA is a cutting-edge approach that harnesses the power of topology, a branch of mathematics that studies the shape and structure of spaces, to analyze and extract valuable insights from data. In the fast-paced evolution of technical achievements, one of the most valuable resources is data. In scientific, economical and political areas, data are being gathered in a progressively fast pace which amounts to an exponential increase of new information each year. Traditional analytical methods often fall short in capturing the underlying patterns and structures buried within these data. This is where TDA steps in with its innovative approach. By considering data as a collection of points in a high-dimensional space, TDA aims to unveil the intrinsic geometry and relationships present in these datasets.

The significance of TDA spans a wide range of fields, making it a transformative tool across diverse disciplines. In data science, TDA offers an alternative approach to dimensionality reduction, helping to projecting datasets onto representations that retain essential information (e.g., see \cite{nicolau2011topology}, \cite{yu2021shape}). In biology, TDA aids in understanding complex biological systems, such as protein structures or neural connectivity, by revealing the underlying topological features that govern their behavior (e.g., \cite{koseki2023topological}, \cite{das2023topological}). Moreover, in image analysis, TDA can decode intricate patterns within images, going beyond pixel-level analysis to uncover hidden structures that might represent critical information (e.g., see \cite{ver2023primer}). This is particularly relevant in medical imaging, where TDA's ability to highlight significant features can aid in disease diagnosis and treatment planning (e.g., \cite{singh2023topological}).

TDA's strength lies in its ability to capture essential characteristics of data that might be overlooked by traditional techniques. By focusing on the shape, connectivity, and arrangement of data points, TDA can identify clusters, holes, voids, and other topological features that convey crucial insights about the data's underlying structure. This capability is especially powerful when dealing with noisy, high-dimensional, or incomplete datasets, where conventional methods often struggle.

The given data types can range from a point cloud to time-series to even visual images and are highly dependent on the subject and method of obtaining said data. In this paper, we focus on spectral data which are typical in fields such as biology and in particular medicine. The most obvious choice for analyzing this type of data utilizing tools of TDA is applying a levelset filtration (e.g., see \cite{edelsbrunner2008persistent}) which extracts important topological information of the underlying space in the form of persistent homology. The results can be represented as barcodes (e.g., see \cite{ghrist2008barcodes}) or persistence diagrams (e.g., see \cite{cohen2005stability}) which can be interpreted individually dependent on the application. However, the levelset filtration has difficulties in differentiating symmetric spectra, sharing the same significant peaks at different locations. In such scenarios, the generated persistence diagrams might appear identical, disrupting the analysis. 

One example for this disturbance can be observed in MALDI imaging, where the importance of positional information in spectral data is linked with the persistence of signal peaks. This insight inspired the closely connected paper "Supervised topological data analysis for MALDI mass spectrometry imaging applications" by G. Klaila, V. Vutov and A. Stefanou (\cite{klaila2023supervised}). In this paper, MALDI imaging was used to classify cancer cells utilizing machine learning and feature extraction with the persistence transformation. This tool helped in pre-processing the signal peaks and their positional information in order to improve the training of the machine learning algorithm. In this example, the position of the peaks could be related to specific molecules which in turn gave information about the type of cancer cell.

Central to our research is the profound significance of stability. While our paper indeed delves into the discrimination of persistence and peaks, a core facet also centers on substantiating the stability of our innovative methodology. In the realm of analytical approaches, stability emerges as a pivotal element, indicating the reliability of outcomes in the face of uncertainties and fluctuations. Stability of the analysis implies robust results related to small perturbations in the input. By demonstrating the stability of our method, we reinforce its credibility and applicability, rendering it a robust tool capable of withstanding the challenges of real-world data.

The first TDA stability theorem by Cohen-Steiner et al. (\cite{cohen2005stability}) proved the stability of persistence diagrams for certain continuous functions. This theorem was expanded to the application of the interleaving distance on persistence modules by Chazal et al. (\cite{chazal2009stability}), and in succession to other applications. Examples for these are the generalizations of the interleaving distance and thus generalizations of the original stability result (see \cite{bubenik2015metrics,stefanou2018interleaving}) or the stability of the Euler characteristic curves and its multi-parameter extension, the Euler characteristic profile (see \cite{dlotko2022euler}). In this paper, we want to contribute to this already established framework of TDA by stating a stability result for the case of the persistence transformation. 

Providing this stability is not just an academic pursuit, but a practical necessity. It enhances the credibility and applicability of our approach, making it a robust tool for real-world data analysis. Our exploration into stability promises to be an essential step forward in the dynamic landscape of modern data analysis, where data complexities are met with steadfast methodologies.

\subsection{Overview of our results}
During the course of this paper, we explore the utilization of the persistence transformation, a novel method in topological data analysis, and delve into its potential benefits. The paper is organized into several chapters, each contributing to a comprehensive understanding of this methodology.

Section \ref{section_morse_sets} provides the necessary background and lays the foundation for our analysis. We define the space of critical points, called Morse Sets, in \ref{subsection_morse_sets} and establish a $\ell_p$-metric to support our subsequent developments in \ref{subsection_lp_metrics}. 

In Section \ref{section_persistence_transformation}, we use this groundwork in order to introduce the persistence transformation in \ref{subsection_persistence_transformation}, showcasing its formulation through the introduction of a matching mechanism on the set of local maxima. Moreover, we place significant emphasis on proving a stability theorem in \ref{subsection_stability}, which stands as the pivotal result of our work. To provide context, we conclude this chapter by providing a brief comparison between the persistence transformation and the conventional persistence diagram in \ref{subsection_comparison}, along with a description of a potential implementation in \ref{subsection_implementation}.

Our exploration continues in Section \ref{section_reduced_persistence_transformation}, where we unveil the reduced persistence transformation - a modified version that offers dimensionality reduction. We carefully describe its structure in \ref{subsection_reduced_definition} and demonstrate its stability in \ref{subsection_reduced_stability}. Additionally, we compare its performance against the traditional persistence diagram in \ref{subsection_reduced_comparison}, shedding light on its strengths and limitations.

Intriguingly, the Section \ref{section_higher_dimensions} offers a conceptual sketch for extending the persistence transformation to higher-dimensional input. We explore the exciting possibilities and inherent challenges that arise in this endeavor, and we discuss the stability of this extension in light of these new conditions.

As we approach the conclusion, in Section \ref{section_summary} we summarize the key findings presented throughout the paper. We emphasize the significance of the persistence transformation as a powerful and stable tool in topological data analysis. Additionally, we give an outlook by considering potential enhancements and novel applications for this transformative method.

Through these chapters, we aim to establish the persistence transformation as a crucial addition to the arsenal of topological data analysis techniques. We showcase its potential to unlock new insights in various domains and illuminate its stability in diverse scenarios.

\section{Morse Sets}\label{section_morse_sets}
\subsection{Background}\label{subsection_background}
In numerous modern applications, significant data are represented as images of real-valued functions (e.g., science, politics, economy, healthcare, engineering), known as time series. These fields include also domains such as medicine (e.g., \cite{radtke2016multiparametric}), robotics (e.g., \cite{atkeson1986robot}), and climate science (e.g., \cite{ge2010temperature}). While studying these functions, one can conclude that relevant information is often given by the critical values of the functions, i.e., the minima and maxima. For example, in the application of MALDI-Imaging (see \cite{aichler2015maldi}), the local maxima of MALDI-Spectra can be associated to specific molecules. This association can help the tumor sub-typing process (see \cite{klaila2023supervised}). On the other hand, values in between these critical points often carry less to none significant information for the analysis of the data. Expanding upon this observation, we establish an abstract space of critical points derived from real-valued functions to facilitate data analysis simplification. In the described context, the abstract space contains similar significant information in a more compact form compared to the original data.
Another notable benefit of this abstraction is that the analysis of the data can be performed without any prior knowledge of the original function, thus enabling the application of these methods to any (data-) set contained within this abstract space. We call these abstract sets \textit{Morse sets}.

\subsection{Definition of Morse sets}\label{subsection_morse_sets}
The term 'Morse set' draws its inspiration from its application in the context of Morse functions. In numerous scenarios, time series data lends itself to interpretation as a real-valued function that satisfies all the prerequisites of a Morse function, as detailed in \cite{milnor1963morse}. However, it is essential to note that the concept of Morse sets can also extend beyond Morse functions and be defined for a broader class of real-valued functions. In this subsection, we provide a formal definition of the Morse set.

For a compact subset, $M\subseteq \overline{\mathbb{R}}:=\mathbb{R}\cup \{-\infty, +\infty\}$, we define the space $\mathscr{M}:=M\times \overline{\mathbb{R}}$ and equip it with the co-lexicographic order $<_{\mathscr{M}}$, which is induced by $<_M$ as follows: 

For all $(x, y), (x', y)$ and $(x', y')\in \mathscr{M}$ with $x\neq x'$ and $y\neq y'$ we have: \label{def_M}
\begin{itemize}
    \item $(x,y)<_{\mathscr{M}} (x', y')\Leftrightarrow y<_{\overline{\mathbb{R}}}y'$
    \item $(x, y)<_{\mathscr{M}} (x', y)\phantom{'}\Leftrightarrow x<_Mx'$.
\end{itemize}
Sets $K$ of elements of this space are denoted by as either $k\in K$ or $(x,y)\in K$, where $x\in M$ and $y\in \overline{\mathbb{R}}$. The set $K\subseteq \mathscr{M}$ is referred to as a \textit{Morse set}, if the following conditions are satisfied:
\begin{enumerate}\label{def_K}
    \item Injectivity: For all $x\in M, y, y'\in \overline{\mathbb{R}}$ there is: $(x, y), (x, y')\in K\Rightarrow y=y'$.
    \item Disjunction: The set $K$ is a disjoint union of two subsets, i.e., $K=K^+\sqcup K^-$. To highlight the affiliation to one of these set, we denote the elements with the corresponding symbol, e.g., $k_i^+\in K^+$ or $(x^-, y^-)\in K^-$.
    \item Ordered: For all $i<j$, it holds that $k_i^+>_{\mathscr{M}}k_j^+$ with $k_i^+, k_j^+\in K^+$ and $k_i^-<_{\mathscr{M}}k_j^-$ with $k_i^-, k_j^-\in K^-$.
    \item Alternation: For all $(x,y),(x', y')\in K^+$ holds: If there is no element $(x^*, y^*)\in K^+$ such that $x^*\in [x,x']$, then there is exactly one element $(\hat{x}, \hat{y})\in K^-$ such that $(\hat{x}, \hat{y})< (x,y),(x', y')$ and $\hat{x}\in [x,x']$, and vice versa.
    
    For all $(x,y),(x', y')\in K^-$ holds: If there is no element $(x^*, y^*)\in K^-$ such that $x^*\in [x,x']$, then there is exactly one element $(\hat{x}, \hat{y})\in K^+$ such that $(\hat{x}, \hat{y})> (x,y),(x', y')$ and $\hat{x}\in [x,x']$.
    \item Critical Boundary: For all $x\in \partial M $ (i.e., the boundary of $M\subset \mathbb{R}$) there is $k=(x,y)\in K$, with $y\in \overline{\mathbb{R}}$.
\end{enumerate}

Let $\mathscr{K}:=\{K\subset \mathscr{M}|K \text{ Morse set}\}$ be the space of all Morse sets. One can think of a Morse set as being the set of non-degenerate critical points for the graph of a Morse function (see \cite{milnor1963morse}) $f:M\rightarrow \overline{\mathbb{R}}$ together with the boundary points, where the $(x, f(x))\in K^+$ are the maxima, and the $(x', f(x'))\in K^-$ are the minima. 

We denote the cardinality of the subsets $K^+, K^-\subset K$ with $\kappa^+:=|K^+|$ and $\kappa^-:=|K^-|$. The condition (4) of \ref{def_K} directly implies the following  statement:
\begin{align*}\label{balanc}
    &&||K^+|-|K^-||&\leq 1\\
    \Leftrightarrow && |\kappa^+ -\kappa^-|&\leq 1.
\end{align*}

A Morse set can be directly obtained by a Morse function $f$ by classifying the non-degenerate critical points in the corresponding subsets $K^+$ and $K^-$. We denote the Morse set $K_f$ to indicate the relationship with $f$. The benefits of employing Morse sets over the original Morse functions become more apparent when considering isotopy, as Morse sets encapsulate information regarding isotopy classes of Morse functions, making them more versatile. Furthermore, when Morse sets are identical, they exhibit shared characteristics related to persistent homology.

According to Kudryavtseva (see \cite{kudryavtseva2009uniform}), two Morse functions are defined to be isotopic, if
\begin{itemize}
    \item there are diffeomorphisms $h_1$ and $h_2$ such that $f = h_2 \circ g \circ h_1$,
    \item  $h_1$ preserves the numbering of critical points
\item $h_1$ is homotopic to the identity mapping.
\end{itemize}

Furthermore, if the critical points are identical, we call the functions \textit{Morse isotopic}. This directly implies the identity of the Morse sets, i.e., if two Morse functions $f$ and $g$ are Morse isotopic, there is $K_f = K_g$. The subsequent theorem demonstrates that working with Morse sets not only encodes more information but also serves as a foundation for deriving concepts such as isotopy and persistent homology for similar Morse functions.

\begin{theorem}
    If two Morse functions $f$ and $g$ are Morse isotopic with $K_f = K_g$, then the following two statements are satisfied:
    \begin{enumerate}
        \item The Morse functions $f$ and $g$ are isotopic.
        \item The persistent homology classes of the upper levelset filtrations are identical.
    \end{enumerate}
\end{theorem}
\begin{proof}
    \begin{enumerate}
    \item For $K_f$, let $G_f$ be the graph obtained by connecting all points $k\in K_f$ with straight lines. The graph $G_f$ is isotopic to the graph $G_g^*$ of $f$, since the straight lines of $G_f$ can continuous be deformed to the graph $G_f^*$. Since the critical points of $G_f$ and $G_f^*$ are identical, the deformation will not bring any more critical points. Similarly, a graph $G_g$ for $K_g$ can be constructed, which is isotopic to the graph $G_g^*$ of $g$. Since $K_f=K_g$, there is equality between $G_f$ and $G_g$, hence there is an isotopy between $G_f^*$ and $G_g^*$. Finally, Theorem 2 of Kudryavtseva (see \cite{kudryavtseva2009uniform}) implies the isotopy of $f$ and $g$.
        \item The persistent homology class of a levelset filtration depends on the path-connected components of the levelsets. For the levelset, a new path-connected component is born, if a critical point is passed. The merging of two path-connected components occurs also during the passing of a critical point. For example. in the upper levelset filtration, a new path-connected component is born in each maxima, and two components are merging in each minima. Since the critical points of $f$ and $g$ are identical, i.e., $K_f=K_g$, the path-connected components of them are also identical in any step. Hence their persistent homology classes are equal, resulting in the same persistent diagram.
    \end{enumerate}
\end{proof}

Utilizing Morse sets allows for the handling of Morse isotopy classes of Morse functions rather than dealing with each Morse function individually. This aids in categorizing Morse functions within a broader context, ensuring that similar Morse functions exhibit comparable behavior in terms of their persistent homology.

\subsection{\texorpdfstring{$\ell_p$}{TEXT}-Metrics on Morse sets}\label{subsection_lp_metrics}
This work aims at proving stability of the persistence transformation, a method working on the above defined space $\mathscr{K}$. Any stability theorem depends on the choice of metrics on the input and output spaces. One metric often utilized when working with functions is the supremum norm, i.e., the maximal distance between two functions, as can be seen in the proof of the stability theorem of the persistence diagram (see \cite{cohen2005stability}). Translated to two sets $K_f, K_g\in \mathscr{K}$, this would be the maximal distance for the minimal matching between elements $k_f\in K_f$ and $k_g\in K_g$. However, the persistence diagram considers stability only in one dimension, therefore it is enough for the metric to control distance in a single dimension. In contrast, the later defined persistence transformation will consider stability not only in the height, but also in the position of the relevant peaks in accordance to the elder rule (see \cite{Edelsbrunner2010}). To satisfy these conditions, we define an $\ell_p$ metric on $\mathscr{K}$ respecting the total order $<_{\mathscr{M}}$ of the elements and controlling the height as well as the position of similar significant peaks in the input. For this metric, these peaks are matched to each other, resulting in a measurable distance for Morse sets.

Define the matching $m^*$ between the Morse sets $K$ and $\hat{K}$ as follows:
\begin{itemize}
    \item For all $i = 1, \ldots, \min\{\kappa^+, \hat{\kappa}^+\}$ there is $(k_i^+, \hat{k}_i^+)\in m^*$.
    \item For all $j = 1, \ldots, \min\{\kappa^-, \hat{\kappa}^-\}$ there is $(k_j^-, \hat{k}_j^-)\in m^*$.
    \item For all $i=\min\{\kappa^+, \hat{\kappa}^+\} + 1, \ldots,\max\{\kappa^+, \hat{\kappa}^+\}$ there is $(0_{\mathscr{M}}, \hat{k}_i^+)\in m^*$ (with $0_{\mathscr{M}}=(0,0)$).
    \item For all $j=\min\{\kappa^-, \hat{\kappa}^-\} + 1, \ldots, \max\{\kappa^-, \hat{\kappa}^-\}$ there is $(0_{\mathscr{M}}, \hat{k}_j^-)\in m^*$. 
\end{itemize}
With this matching $m^*$, an $\ell_p$ metric on $\mathscr{K}$ can be established.
\begin{definition}\label{metric_M}
For all $K, \hat{K}\in \mathscr{K}$ with the matching $m^*$, the distance between $K$ and $\hat{K}$ is defined to be:
\begin{align*}
    d_{\mathscr{K}, p}(K, \hat{K}):&=\sqrt[p]{\sum_{(k, \hat{k})\in m^*}||k-\hat{k}||_{\infty}^p}\\
    &=\sqrt[p]{\sum_{((x,y), (\hat{x},\hat{y}))\in m^*}||(x,y)-(\hat{x},\hat{y})||_{\infty}^p}.
\end{align*}
\end{definition}
\begin{theorem}
    The metric $d_{\mathscr{K}, p}(K, \hat{K})$ is well-defined.
\end{theorem}

\begin{proof}
    The metric $d_{\mathscr{K}, p}(K, \hat{K})$ is well-defined since the following four conditions are satisfied:
    \begin{itemize}
        \item Non-negativity. For all $K, \hat{K}\in \mathscr{K}$, there is:
        \begin{align*}
            d_{\mathscr{K}, p}(K, K')&=\sqrt[p]{\sum_{(k,\hat{k})\in m^*}||k-\hat{k}||_{\infty}^p}\\
            &\geq \sqrt[p]{\sum_{(k, \hat{k})\in m^*}0}\\
            &= \sqrt[p]{0}\\
            &=0,
        \end{align*}
        since $||\cdot , \cdot||_{\infty}$ is a metric and therefore non-negative.
        \item Definiteness. For all $K, \hat{K}\in \mathscr{K}$ there is: 
        \begin{align*}
            K=\hat{K} \Rightarrow d_{\mathscr{K}, p}(K, \hat{K})=0,
        \end{align*}
        since $m^*$ is the identity matching on $K=\hat{K}$. On the other hand, assume that
        \begin{align*}
            d_{\mathscr{K}, p}(K, \hat{K})=0.
        \end{align*}
        Then there is a matching $m^*$, such that for all $i=1, \dots, \min \{\kappa^+, \hat{\kappa}^+\}$  there is $(k_i^+, \hat{k_i}^+)\in m^*$ and for all $j=1, \dots, \min\{\kappa^-, \hat{\kappa}^-\}$ there is $(k_j^-, \hat{k_j}^-)\in m^*$. The zero distance between $K$ and $\hat{K}$ implies equality between all these elements, i.e., $k_i^+=\hat{k_i}^+$ and $k_j^-=\hat{k_j}^-$. Assume now without loss of generality that $\kappa^+\geq \hat{\kappa}^+$ and $\kappa^-\geq \hat{\kappa}^-$. It holds, that for all $i=\hat{\kappa}^++1, \dots, \kappa^+$ the elements $k_i^+\in K^+$ are matched to $0_{\mathscr{M}}$, and similar for all $j=\hat{\kappa}^-+1, \dots, \kappa^-$ the elements $k_j^-\in K^-$ are matched to $0_{\mathscr{M}}$. However, since the distance between $K$ and $\hat{K}$ is zero, the elements must be equal to $0_{\mathscr{M}}$. Uniqueness (\ref{def_K}, (1)) of the sets implies the existence of at most one such element $0_{\mathscr{M}} \in K$. Consequently, there must be equality between all elements $k\in K$ and $\hat{k}\in \hat{K}$, except for at most one element $0_{\mathscr{M}} \in K$. Without loss of generality let $0_{\mathscr{M}}\in K^+$. This element cannot be in the closure of $M$, because the critical boundary (\ref{def_K}, (5)) of $\hat{K}$ implies the existence of $\hat{k}\in \hat{K}$ which is in the boundary. Equality of the elements then states the existence of an corresponding element $k\in K$ in the boundary, being different than $0_{\mathscr{M}}$ and hence contradicting the uniqueness (\ref{def_K}, (1)) of the boundary element. If the element $0_{\mathscr{M}}$ is not in the boundary of $M$, the alternation of $K$ (\ref{def_K}, (4)) implies the existence of neighboring elements $k=(x,y), k'=(x', y')\in K^-$ such that $(0_M, 0_{\mathbb{R}})=0_{\mathscr{M}}\in K^+$ is a unique element in $K^+$ with $0_M\in [x,x']$. Per assumption, the elements of $K$ and $\hat{K}$ are the same except for $0_{\mathscr{M}}$, so there are elements $\hat{k}, \hat{k}'\in \hat{K}^-$ such that $k=\hat{k}$ and $k'=\hat{k}'$. The alternation of $\hat{K}$ (\ref{def_K}, (4)) implies the existence of at least one element $\hat{k}^*=(\hat{x}^*, \hat{y}^*)\in \hat{K}^+$ such that $\hat{x}^*\in [x, x']$. By equality of the elements, there must be an element $k^*=(x^*, y^*)\in K^+$ such that $k^*\neq 0_{\mathscr{M}}$ and $x^*\in [x, x']$, opposing the uniqueness of $0_{\mathscr{M}}$ in this interval. This contradicts the assumption, that there is an element $0_{\mathscr{M}}\in K$ with $0_{\mathscr{M}}\notin \hat{K}$, hence all the elements must be the same, proofing:
        \begin{align*}
             d_{\mathscr{K}, p}(K, \hat{K})=0 \Rightarrow K=\hat{K}.
        \end{align*}
        \item Symmetry. For all $K, \hat{K}\in \mathscr{K}$ there is: 
        \begin{align*}
            d_{\mathscr{K}, p}(K, \hat{K})^p&=\sum_{(k, \hat{k})\in m^*}||k-\hat{k}||_{\infty}^p\\
            &=\sum_{(k, \hat{k})\in m^*}||\hat{k}-k||_{\infty}^p\\
            &=\sum_{(\hat{k}, k)\in m^*}||\hat{k}-k||_{\infty}^p\\
            &= d_{\mathscr{K}, p}(\hat{K}, K)^p.
        \end{align*}
        \item Triangle inequality. For all $K, \hat{K}, \Tilde{K}\in \mathscr{K}$ assume without loss of generality $k^+=\hat{k}^+=\Tilde{k}^+$ and $k^-=\hat{k}^-=\Tilde{k}^-$. If not, add elements $0_{\mathscr{M}}$ to the sets with fewer elements. Then:
        \begin{align*}
            &d_{\mathscr{K}, p}(K, \Tilde{K})^p\\
            =&\sum_{(k, \Tilde{k})\in m^*_{K, \Tilde{K}}}||k-\Tilde{k}||_{\infty}^p\\
            =&\sum_{i =1}^{\kappa^+}||k^+_i-\Tilde{k}^+_i||_{\infty}^p+\sum_{j=1}^{\kappa^-}||k_j^--\Tilde{k}^-_j||_{\infty}^p\\
            = &\sum_{i =1}^{\kappa^+}||k^+_i-\hat{k}^+_i+\hat{k}^+_i -\Tilde{k}^+_i||_{\infty}^p+\sum_{j=1}^{\kappa^-}||k_j^--\hat{k}^-_i+\hat{k}^-_i-\Tilde{k}^-_j||_{\infty}^p\\
            \leq &\sum_{i =1}^{\kappa^+}||k^+_i-\hat{k}^+_i||_{\infty}^p+||\hat{k}^+_i -\Tilde{k}^+_i||_{\infty}^p+\sum_{j=1}^{\kappa^-}||k_j^--\hat{k}^-_i||_{\infty}^p+||\hat{k}^-_i-\Tilde{k}^-_j||_{\infty}^p\\
            =& \sum_{(k, \hat{k})\in m^*_{K, \hat{K}}}||k-\hat{k}||_{\infty}^p+\sum_{(\hat{k}, \Tilde{k})\in m^*_{\hat{K}, \Tilde{K}}}||\hat{k}-\Tilde{k}||_{\infty}^p\\
            =& d_{\mathscr{K}, p}(K, \hat{K})^p+d_{\mathscr{K}, p}(\hat{K}, \Tilde{K})^p.
        \end{align*}
    \end{itemize}
\end{proof}
\section{Persistence Transformation}\label{section_persistence_transformation}
One main objective of topological data analysis is to track the persistence of topological features in the data. Each feature gets assigned a birth value, i.e., the value where the feature arises, and a death value, i.e., the value where the feature merges to other features. In the persistence transformation, the merging process is according to the elder rule (see \cite{Edelsbrunner2010}). The persistence of a feature is the lifespan of it, in other words, the difference of birth and death. 

\subsection{Matching}\label{subsection_matching}
In this paper, we consider features to be the peaks of a Morse function, or in the notation from before, the elements $k^+=(x^+, y^+)\in K^+$, i.e., the local maxima. The birth of these features is given natural by their height, i.e., by $y^+$. On the other hand, the merging of two features always happens at a local minimum, or the $(x^-, y^-)=k^-\in K^-$. For each feature, there is a unique $k^-\in K^-$, where it dies, and the death value is given by $y^-$. To track the death of the feature, we match to each feature $k^+\in K^+$ the unique element of merging $k^-\in K^-$. The largest feature cannot be merged with another feature according to the elder rule (see \cite{Edelsbrunner2010}), therefore its death is defined to be $-\infty$. The matching is denoted by $\mu:K^+\rightarrow K^-\cup M\times \{-\infty\}$ and is defined sequential and individual on each path-connected component. The results can later be joined for the persistence transformation.

\begin{itemize}\label{def_mu}
    \item For all $i=\kappa^+, \dots, 2$:
    \begin{align*}
        \mu(k_i):=\sup &\{(x^-, y^-)=k^-\in K^-|k^-< k_i\\
        \land &\exists k_i\leq k^+=(x^+,y^+)\in K^+: x^-\in [\min\{x_i, x^+\}, \max \{x_i, x^+\}]\\
        \land & \forall j = i+1, \dots, \kappa^+: k^-\neq \mu(k_j) \}.
    \end{align*}
    \item For $i=1$:
    \begin{align*}
        \mu(k_1):=(x_1, -\infty).
    \end{align*}
\end{itemize}

This functional determination approach presents an innovative method for calculating persistence's in accordance with the elder rule by utilizing a semi-algorithmic procedure, which departs from the conventional complex algorithms typically employed. An important aspect of this is the simplification of the process of obtaining persistence's, making it more accessible and efficient.

\begin{theorem}
\label{injective}
    The matching $\mu$ is injective and is well-defined.
\end{theorem}
\begin{proof}
    Via natural induction over $\kappa^+$. 
    
    Induction start for $\kappa^+=2$, i.e. $K^+=\{k_1, k_2\}$:
    The Alternation (\ref{def_K}, (4)) implies the existence of an element $k^-\in K^-$, such that $k^-< k_2$ and $x^-\in [\min \{x_1, x_2\}, \max \{x_1, x_2\}]$. Hence this element is suitable for $\mu(k_2)$. The element $k_1$ is matched to $\mu(k_1)=(x_1, -\infty)$, resulting in an injective matching.

    For the induction step assume that we can always find an injective matching for $\kappa^+=n$. For $\kappa^+=n+1$ now holds: The Alternation (\ref{def_K}, (4)) implies with the same reasoning of the induction start the existence of an suitable element $k^-$, such that $\mu(k_{\kappa^+})=k^-$. This element is by (\ref{def_K}, (4)) a unique element of $K^-$ for a suitable $k_i\in K$ with $x^-\in [\min \{x_{\kappa^+}, x_i\}, \max \{x_{\kappa^+}, x_i\}]$, such that there are no other elements $k^+\in K^+$ in that interval. Therefore, the points $k^-, k_{\kappa^+}$ can be removed from the set $K$, resulting in a set $\hat{K}$ with $\hat{\kappa}^+=n$, which satisfies all conditions of an ordered critical set (\ref{def_K}). The induction assumption applied to $\hat{K}$ completes the proof.
\end{proof}

\subsection{Definition of the Persistence Transformation}\label{subsection_persistence_transformation}
\begin{figure}
    \centering
    \includegraphics{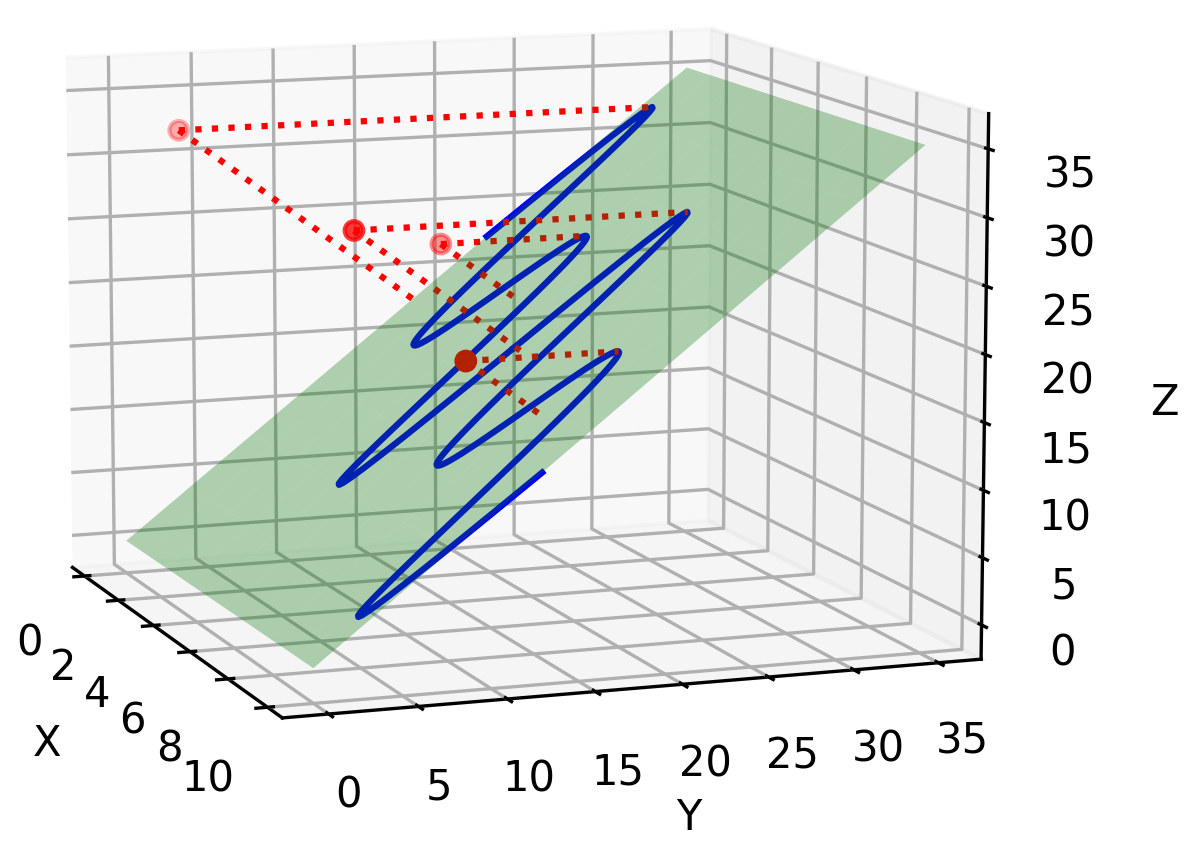}
    \caption{Example of the persistence transformation. The blue line represents all the trivial features, which are vanishing since they are on the diagonal plane. The red dots represent the relevant features. The distance of the dots to the diagonal plane indicates their persistence.}
    \label{fig:example_pt}
\end{figure}

With the defined matching, the persistence transformation can be defined. It is constructed in such a way, that it tracks the position, the birth and the death of each feature $k\in K^+$.
\begin{definition}\label{def_persistence_transformation}
    The \textbf{persistence transformation} is a map
    \begin{align*}
        T:\mathscr{K}&\rightarrow \mathscr{M}\times \overline{\mathbb{R}}\\
        K &\mapsto T_K,
    \end{align*}
    where each element $k=(x,y)\in K^+$ with $\mu(k)=\overline{k}=(\overline{x}, \overline{y})$ is mapped to $t_K:=(x,y,\overline{y})$. The elements $k^-=(x^-, y^-)\in K^-$ are mapped to $t_{K^-}=(x^-, y^-, y^-)$ on a diagonal plane.
\end{definition}
A graphical representation of the persistence transformation can be seen in Fig. \ref{fig:example_pt}. The first coordinate of each element $t_K\in T_K$ corresponds to the position $x\in M$ of a feature. The second coordinate denotes the birth value $y\in\overline{\mathbb{R}}$, while the last coordinate describes the death value, i.e., $\overline{y}\in \overline{\mathbb{R}}$. This results in trivial features for all elements $k^-\in K^-$, since their death value equals their birth value. They are located on the diagonal plane in Fig. \ref{fig:example_pt}, which can be compared to the diagonal in the persistence diagram.

\subsection{Stability}\label{subsection_stability}
A necessary characteristic of a useful method in topological data analysis is the stability of the method. The stability is a measure of how changes in the input influence the output of the analysis and indicates in this sense, how robust the applied method captures the topological information of the data. A stable method ensures that the variation in the output is bound by the variation in the input (see \cite{cohen2005stability}). The subsequent section provides evidence for the stability of the persistence transformation. For this, the $p$-Wasserstein distance is utilized a the metric on the output of the persistence transformation. This distance is defined as follows (see \cite{memoli2011gromov, berwald2018computing}): 
 \begin{align*}
            d_{W_p}(A, B)^p =\min_{\text{matchings }m:A\times B}\sum_{(a,b)\in m}||a-b||_{\infty}^p.
\end{align*}
As $p$ approaches infinity, this metric converges to the bottleneck distance, the metric employed in assessing the stability of persistence diagrams (refer to \cite{cohen2005stability}). The bottleneck distance exhibits notable computational efficiency and robustness to outliers when compared to finite $p$-Wasserstein distances. Furthermore, it offers a more intuitive interpretation by quantifying the maximum separation between elements. Importantly, by establishing the theorem's validity for arbitrary values of $p$, we ensure that the obtained results are equally applicable to the bottleneck distance.

\begin{theorem}
    The persistence transformation is stable using the $p$-Wasserstein distance on $\mathscr{M}\times \overline{\mathbb{R}}$:
    \begin{align*}
        d_{W_p}(T_K, T_L)\leq d_{\mathscr{K}, p}(K, L).
    \end{align*}
\end{theorem}
\begin{proof}
Let $P=T_K\times T_L$ be the set of all matchings for $t_K=(x,y,\hat{y})\in T_K$ and $t_L=(x', y', \overline{y'})\in T_L$, given $k=(x,y)\in K$ and $l=(x', y')\in L$. Without loss of generality let $\kappa^+:=\max\{\kappa_K^+, \kappa_L^+\}$ and $\kappa^-=\max\{\kappa_K^-,\kappa_L^-\}$. The stability of the persistence transformation is given by the following inequality:
    \begin{align}
        &d_{W_p}(T_K, T_L)^p\\
        =&\min_{m\in P}\sum_{(t_K,t_L)\in m}||t_K-t_L||_{\infty}^p\\
        =&\min_{m\in P}\sum_{(t_K, t_L)\in m}||(x,y, \overline{y})-(x', y', \overline{y'})||_{\infty}^p\\
        =&\min_{m\in P}\sum_{(t_K, t_L)\in m}\max\{||(x,y)-(x', y')||_{\infty}^p, ||\overline{y}-\overline{y'}||_{\infty}^p\}\\
        \leq &\min_{m\in P}\sum_{(t_K, t_L)\in m}\max\{||(x,y)-(x', y')||_{\infty}^p, ||(\overline{x}, \overline{y})-(\overline{x'},\overline{y'})||_{\infty}^p\}\\
        \leq &\sum_{(k, l)\in m^*}\max\{||(x,y)-(x', y')||_{\infty}^p, ||(\overline{x}, \overline{y})-(\overline{x'},\overline{y'})||_{\infty}^p\}\\
        \overset{\ref{injective}}{\leq} &\sum_{i=1}^{\kappa^+}||k_i^+-l_i^+||_{\infty}^p +\sum_{j=1}^{\kappa^-} ||k_j^--l_j^-||_{\infty}^p\\
        =&\sum_{(k,l)\in m^*}||k-l||_{\infty}^p\\
        =& d_{\mathscr{K}, p}(K, L)^p.
    \end{align} 
    The inequality of (7) holds due the injectivity (\ref{injective}) of the matching function. This implies, that each element $k^-\in K^-$ and each element $l^-\in L^-$ are matched at most once.
\end{proof}

\subsection{Comparison}\label{subsection_comparison}

\begin{figure}[ht]
    \centering
    \includegraphics[scale = .27]{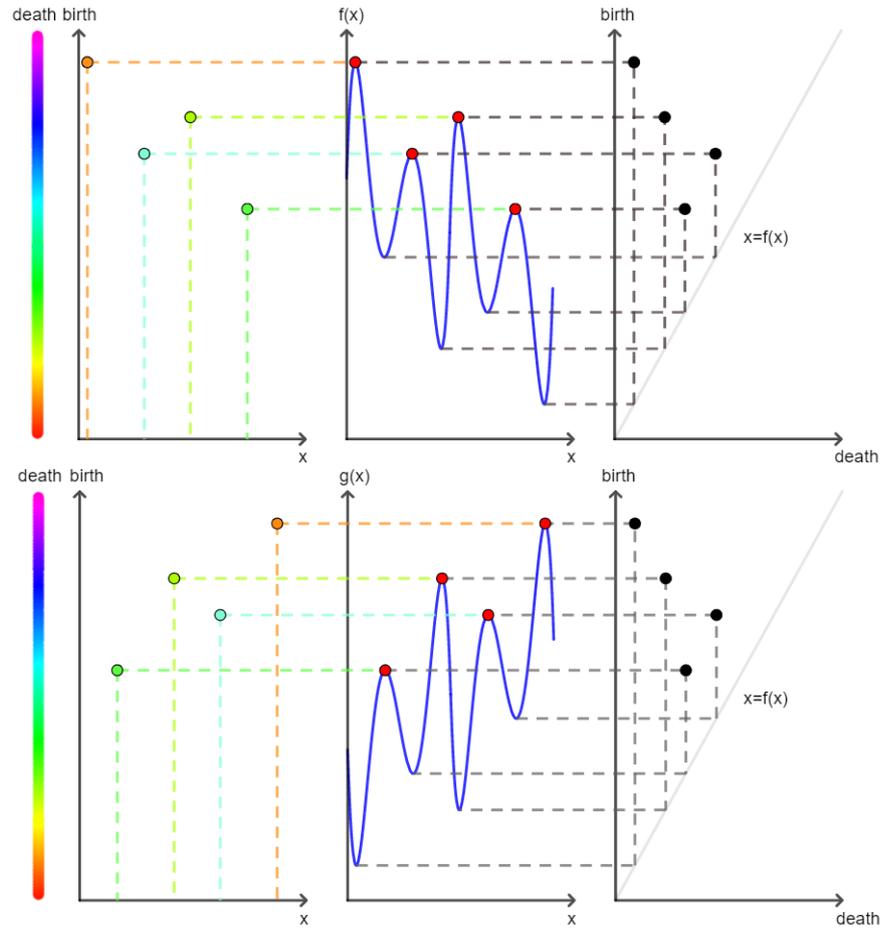}
    \caption{Problem of the persistence diagram: the middle column displays two symmetric functions $f,g:\mathbb{R}\rightarrow \mathbb{R}$. The first column displays the persistence transformation of these functions, with the positional information on the $x$-axis and the birth values on the $y$-axis. The death values are encoded in the color scheme. The persistence transformations of the functions can be distinguished. The right-most column shows the persistence diagram of the upper levelset filtration with the birth value on the $y$-axis and the death value on the $x$-axis. The two diagrams are identical. (Graphic: \cite{klaila2023supervised})}
    \label{fig:xample_pd_fail}
\end{figure}

We conclude this section with a comprehensive comparison between the introduced persistence transformation and the persistence diagram. The persistence diagram captures the information of the persistence in a compact and comprehensive way by displaying the homology classes of specific filtration complexes, e.g., the sub levelset filtration or the upper levelset filtration of a real-valued function. For more information about the standard persistence diagram we refer the reader to \cite{cohen2005stability, edelsbrunner2002topological}. The method has demonstrated successful applications in various occasions (e.g., \cite{edelsbrunner2002topological, otter2017roadmap, nicolau2011topology}). Nonetheless, in some settings, e.g., MALDI-Images (\cite{klaila2023supervised}), it exhibits a significant drawback: while tracking the persistence of the signal peaks, the persistence diagram forfeits the information of the position of the signal. For instance, the standard persistence diagram of the upper levelset filtration cannot distinguish symmetric inputs, see Figure \ref{fig:xample_pd_fail}.

In contrast to the persistence diagram of the upper levelset filtration, the persistence transformation is specifically designed to overcome these limitations of the persistence diagram by tracking the position of significant signal peaks as well as their persistence. Even more, the persistence transformation captures all the information the persistence diagram of the upper levelset filtration captures. This implies that it is a strictly stronger method of capturing the necessary information for distinguishing data. However, this increase in performance comes with the cost of higher dimensionality. While the results of the persistence diagram are subsets of $\overline{\mathbb{R}}^2$, i.e., information about the birth and the death, the results of the persistence transformation are subsets of $M\times \overline{\mathbb{R}}^2$. This poses a significant challenge to the desire of improving computational efficiency for analysis. A decision must be made regarding whether to prioritize fast computation or more comprehensive information. In the upcoming section, we introduce a modification to the persistence transformation that decreases its dimensionality, albeit at the expense of some information.

\subsection{Implementation}\label{subsection_implementation}
In this section, we present the pseudo-code for a potential implementation of the persistence transformation. Given the critical points $K$ as input, the algorithm efficiently computes the output $T_K$ with quadratic complexity. The resulting $T_K$ can then be further processed in linear complexity to represent the reduced persistence transformation of the next section or the persistence diagram. Additionally, the output is sorted based on persistence values. By considering low persistence features as noise, it becomes feasible to denoise the output by removing elements below a specific threshold.

For a similar algorithm, its proof and its practical application, refer to the work in (\cite{klaila2023supervised}).

\begin{algorithm}
\caption{Recursion start}\label{algorithm:recursion_start}
\begin{algorithmic}[1]
\STATE \textbf{Input:} $K^+=\{k_0^+, \ldots, k_{n}^+\}; K^-=\{k_0^-, \ldots, k_{n\pm 1}^-\}$
\STATE \textbf{Return:} $T_K=\{(x_i,y_i,\overline{y_i})|i=0,\ldots,n \}$ 
\STATE $T_K\leftarrow \varnothing$ 
\STATE $k'=(x', y') \leftarrow K^+\text{.pop}(0)$
\STATE $T_K\leftarrow (x', y',-\infty)$
\STATE setOne, setTwo $\leftarrow \varnothing$ 
\STATE minimum $\leftarrow \infty$
\STATE maximum $\leftarrow -\infty$ 
\FORALL{$k=(x,y)\in K^+$}
    \IF{$x<x'$}
        \STATE setOne $\leftarrow k$
        \IF{$x<$ minimum}
            \STATE minimum $\leftarrow x$
        \ENDIF
    \ELSE
        \STATE setTwo $\leftarrow k$
        \IF{$x>$ maximum}
            \STATE maximum $\leftarrow x$
        \ENDIF
    \ENDIF
\ENDFOR
\STATE RecursionStep(minimum$, x'$, setOne, $K^-$.copy(), $T_K$)
\STATE RecursionStep(maximum$, x'$, setTwo, $K^-$.copy(), $T_K$)
\RETURN $T_K$
\end{algorithmic}
\end{algorithm}

\begin{algorithm}
\caption{Recursion Step}\label{algorithm:recursion_step}
\begin{algorithmic}[1]
\STATE \textbf{Input:} start, end, $K^+$, $K^-$, $T_K$
\FORALL{$(x,y)=k\in  K^+$ }
    \IF{$x\notin [\text{start}, \text{end}]$}
        \STATE $K^+\leftarrow K^+\backslash k$
    \ENDIF
\ENDFOR
\IF{$|K^+|=0$}
    \RETURN
\ENDIF
\STATE $k'=(x', y')\leftarrow K^+\text{.pop}(0)$
\STATE RecursionStep(start, $x'$, $K^+$.copy(), $K^-$.copy(), $T_K$)
\FORALL{$(x, y)\in K^-$}
    \IF{$x\notin (x, \text{end})$}
        \STATE $K^-\leftarrow K^-\backslash k$
    \ENDIF
\ENDFOR
\STATE $k^-=(x^-,y^-)\leftarrow K^-$.pop$(0)$
\STATE $T_K\leftarrow (x', y', y^-)$
\STATE RecursionStep($x^-, x'$, $K^+$.copy(), $K^-$.copy(), $T_K$)
\STATE RecursionStep($x^-$, end, $K^+$.copy(), $K^-$.copy(), $T_K$)
\end{algorithmic}
\end{algorithm}

\section{Reduced Persistence Transformation}\label{section_reduced_persistence_transformation} 
\subsection{Motivation and Definition}\label{subsection_reduced_definition}
As stated in \ref{subsection_comparison}, the advantages of the persistence transformation are accompanied by an increase in dimensionality. However, in certain applications, dealing with large amounts of data, the subsequent higher dimensionality in the analyzed data results in prolonged computational time, despite not requiring the full extent of the more comprehensive information gathered. Regardless, these application could also benefit from the positional information provided by the persistence transformation. In response to this concern, we develop an adapted version of the persistence transformation known as the reduced persistence transformation. This modified method aims at reducing the dimensionality of the output data while maintaining the positional information of signal peaks. To accomplish the intended adaptation, we opt for storing the persistence of the features, given by the difference of birth and death, instead of storing these values separately. 

\begin{definition}\label{def_reduced_persistence_transformation}
    Building upon the previous definitions of $\mathscr{M}$ (\ref{def_M}) and the matching $\mu$ (\ref{def_mu}), we define the \textbf{reduced persistence transformation} as a map 
    \begin{align*}
        \Tilde{T}:\mathscr{K}&\rightarrow \mathscr{M}\\
        K^+&\mapsto \Tilde{T}_K
    \end{align*}
    where each element $k=(x,y)\in K^+$ with $\mu(k)=\overline{k}=(\overline{x}, \overline{y})$ is mapped to $\Tilde{t}_K:=(x,y-\overline{y})$.
\end{definition}

\begin{figure}
    \centering
    \includegraphics[scale = .44]{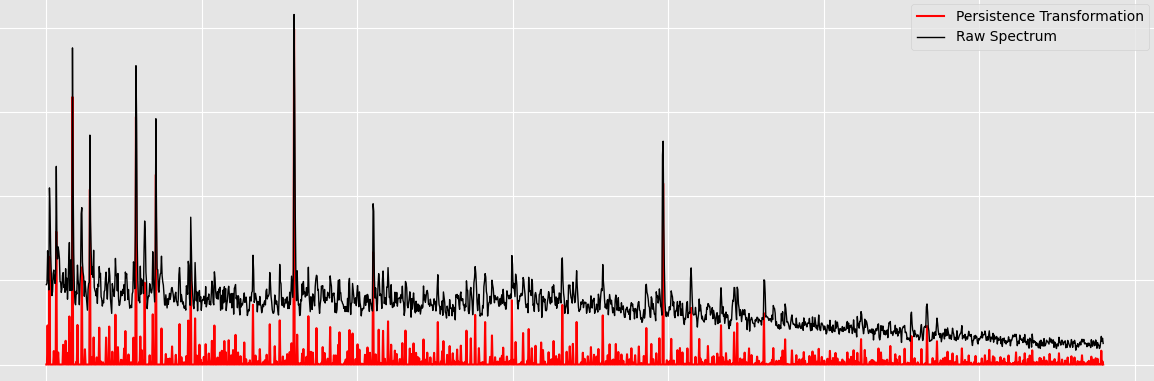}
    \caption{Example of the reduced persistence transformation. The original spectra is displayed in black, while the values of the reduced persistence transformation are illustrating in red the persistence of each feature.}
    \label{fig:example_rpt}
\end{figure}

An example of the reduced persistence transformation can be seen in \ref{fig:example_rpt}

\subsection{Stability}\label{subsection_reduced_stability}
Similar to the stability of the persistence transformation in \ref{subsection_stability}, the stability theorem can also be applied to the reduced persistence transformation.
\begin{theorem}
    The reduced persistence transformation satisfies the stability condition using the $p$-Wasserstein distance in $\mathscr{M}$:
    \begin{align*}
        d_{W_p}(\Tilde{T}_K,\Tilde{T}_L)\leq d_{\mathscr{K}, p}(K, L).
    \end{align*}
\end{theorem}
\begin{proof}
   Let $P=\Tilde{T}_K\times \Tilde{T}_L$ be the set of all matchings for $\Tilde{t}_K=(x,y-\hat{y})\in \Tilde{T}_K$ and $\Tilde{t}_L=(x', y'- \overline{y'})\in \Tilde{T}_L$, given $k=(x,y)\in K$ and $l=(x', y')\in L$. Without loss of generality, let $\kappa^+:=\max\{\kappa_K^+,\kappa_L^+\}$ and $\kappa^-:=\max\{\kappa_K^-,\kappa_L^-\}$. The stability of the reduced persistence transformation is given by the following inequality:
    \begin{align}
        &d_{W_p}(\Tilde{T}_K, \Tilde{T}_L)^p\\
        =&\min_{m\in P}\sum_{(\Tilde{t}_K,\Tilde{t}_L)\in m}||\Tilde{t}_K-\Tilde{t}_L||_{\infty}^p\\
        =&\min_{m\in P}\sum_{(\Tilde{t}_K,\Tilde{t}_L)\in m}||(x,y- \overline{y})-(x', y'- \overline{y'})||_{\infty}^p\\
        \leq & \min_{m\in P}\sum_{(\Tilde{t}_K,\Tilde{t}_L)\in m}||(x,y)-(x', y')||_{\infty}^p+||\overline{y}-\overline{y'}||^P_{\infty}\\
        \leq & \min_{m\in P}\sum_{(\Tilde{t}_K,\Tilde{t}_L)\in m}||(x,y)-(x', y')||_{\infty}^p+||(\overline{x}, \overline{y})-(\overline{x'}, \overline{y'})||^P_{\infty}\\
        \leq & \sum_{(k,l)\in m^*}||(x,y)-(x', y')||_{\infty}^p+||(\overline{x}, \overline{y})-(\overline{x'}, \overline{y'})||^P_{\infty}\\
         \overset{\ref{injective}}{\leq} &\sum_{i=1}^{\kappa^+}||k_i^+-l_i^+||_{\infty}^p +\sum_{j=1}^{\kappa^-} ||k_j^--l_j^-||_{\infty}^p\\
        =&\sum_{(k,\hat{k})\in m^*}||k-l||_{\infty}^p\\
        =& d_{\mathscr{K}, p}(K, L)^p.
    \end{align} 
    The inequality of line (16) is given by the injectivity (\ref{injective}) of the matching function, which ensures the unique occurrence of each element $k^-\in K^-$ and $l^-\in L-$.
\end{proof}

\subsection{Comparison}\label{subsection_reduced_comparison}
\begin{figure}
    \centering
    \includegraphics[scale = .27]{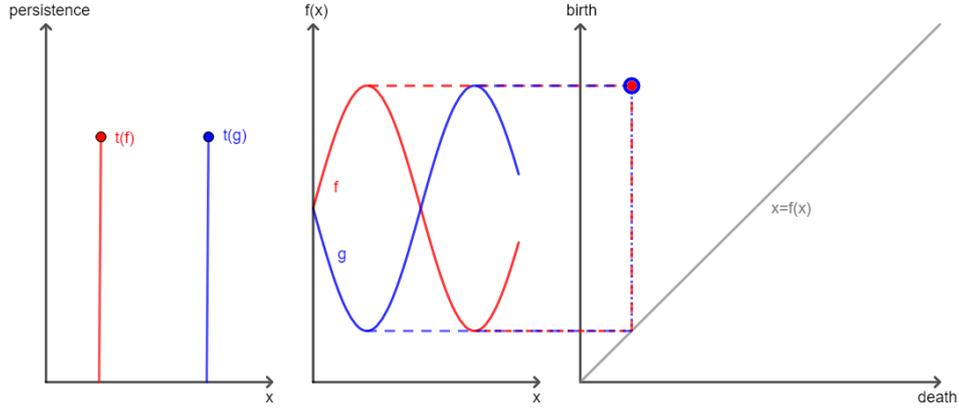}
    \caption{Example of the reduced persistence transformation: In the middle, two symmetric functions $f,g:\mathbb{R}\rightarrow \mathbb{R}$ are displayed. The left-hand side illustrates the reduced persistence transformation of this graphs, with the positional information on the $x$-axis and the persistence information on the $y$-axis. The graphs can be distinguished. The right-hand side depicts the persistence diagram of the upper levelset filtration, with the birth value on the $y$-axis and the death value on the $x$-axis. The diagrams are identical. (Graphic: \cite{klaila2023supervised})}
    \label{fig:example_reduced}
\end{figure}

As mentioned in \ref{subsection_comparison}, while the persistence transformation provides more comprehensive information than the persistence diagram of the upper levelset, it comes at the cost of increased dimensionality. Conversely, the reduced persistence transformation reduces this dimensionality but loses some information in the process. Consequently, it is crucial to conduct a comparison between the reduced persistence transformation and the persistence diagram. After careful examination, it becomes evident that neither of the mentioned methods outperforms the other, and they share the same dimensionality in the output, i.e., $\mathbb{R}^2$. Giving them a total order is not feasible, as there exist scenarios where one outperforms the other, and vice versa (see Fig. (\ref{fig:example_reduced}). The key distinction of the persistence diagram lies in the total height of birth and death values of a feature, whereas the reduced persistence transformation incorporates the positional information alongside its relative height, i.e., its persistence. This implies that in any application where positional information is crucial, the reduced persistence transformation is a more suitable choice as an analysis method. An instance of such an application can be observed in \cite{klaila2023supervised}, where the reduced persistence transformation is employed to enhance the accuracy of tumor tissue classification. In novel applications, it is imperative to assess the significance of positional relevance to determine the potential benefits of employing the modified methodology in comparison to the original persistence transformation or the persistence diagram as each approach possesses its distinct area of efficiency.

\section{Extension to higher dimensions}\label{section_higher_dimensions}
Until this point, we have made the assumption that $M\subseteq \overline{\mathbb{R}}$; however, it is feasible to relax this condition. In this section, we present a concise overview of the approach to extend the persistent transformation to any arbitrary higher-dimensional compact metric space $M$, equipped with a total order and a notion of path connectivity. Additionally, $M$ must be a subset of a space including a neutral element, which is essential for the defined metric on $\mathscr{M}$ (\ref{metric_M}). This requirements can be satisfied by $M\subseteq \overline{\mathbb{R}}^n$, if we equip $\overline{\mathbb{R}}^n$ with a total order, e.g., the lexicographic order. 

One crucial adaptations for the increase of dimension is the definition of the ordered critical set. While in $\overline{\mathbb{R}}$ each path $[x,x']$ is unique for all elements $x,x'\in\overline{\mathbb{R}}$, there may exist multiple paths for elements $x,x'\in \overline{\mathbb{R}}^n$. This property breaks the uniqueness of elements $\hat{k}=(\hat{x},\hat{y})\in K$ such that $\hat{x}\in [x, x']$, resulting in cases where no minimal element $k\in K^-$ lies on the path.

Another issue that needs to be addressed is the presence of areas with the same height, where all elements within that area are critical compared to any element outside that area. Such situations could lead to an imbalance in the number of elements in the sets $K^+$ and $K^-$. To address this, we can utilize the total order of $M$ to select only one element from this area, thus ensuring the satisfaction of the balance equation (\ref{balanc}).

Finally, we also need to consider the adaptation of saddle points, which are degenerate critical points in Morse theory (see \cite{milnor1963morse}). Although they may not possess a unique gradient direction, they can still act as the highest valley that needs to be traversed to reach the next peak, or in other words, they could be crucial for defining the persistence of peaks. One possibility of handling them is to consider these saddle points as maxima and minima, including them in the sets $K^+$ and $K^-$, respectively. In this approach, elements $k \in K$ may have a multiplicity to account for their presence in the analysis.

After implementing all adaptations of the space of critical points, the definitions \ref{metric_M}, \ref{def_mu}, \ref{def_persistence_transformation} and  \ref{def_reduced_persistence_transformation} will hold, resulting in a stable persistence transformation on the higher dimensional set $M$. The proofs of the theorems are presented in a way that they are independent of the dimensionality of $M$. However, for the sake of simplicity and clarity, the specific details of the extension to higher dimensions are omitted in this paper.

It is important to note that increasing the dimension of $M$ will also lead to an increase in the dimension of the output space. As previously mentioned, the reduced persistence transformation and the persistence transformation produce an output of dimensions $M\times \overline{\mathbb{R}}$ and $M\times \overline{\mathbb{R}}^2$, respectively. For cases where $M\subseteq \overline{\mathbb{R}}$, the output dimension of the reduced persistence transformation can match the output dimension of the persistence diagram, i.e., $\overline{\mathbb{R}}^2$. However, for any $M\subseteq \overline{\mathbb{R}}^n$, the output dimension of the reduced persistence transformation will have an increased in dimension of $n-1$, One should consider the potential increase in computational time for further processing steps and carefully evaluate the significance of the positional information obtained from the increased output dimensions.

\section{Summary and Outlook}\label{section_summary}

In conclusion, this paper has successfully achieved its research objectives by constructing a stable method for analyzing time series data arising from Morse functions in a topological manner while preserving the positional information of signal peaks. We have defined the persistence transformation, established its stability, and demonstrated that it complements nicely the existing machinery of persistent homology on levelset filtrations. Furthermore, we have introduced the reduced persistence transformation as a valuable side result, effectively tracking positional information with reduced dimensionality.

The implications of this research are far-reaching, as the proposed method can be applied across various fields where data is represented as the image of a real-valued function. Its potential to deliver promising results, particularly in domains where positional information is crucial, underscores its significance.

Our hypothesis regarding the successful application of the reduced persistence transformation to real-world problems has been confirmed in a previous study about MALDI-MSI data in (\cite{klaila2023supervised}).

However, we acknowledge certain limitations in the persistence transformation, particularly in terms of its computational complexity. As the dimension of the data-set M increase, the output dimension of the persistence transformation also grows, resulting in longer computational times for subsequent processing tasks.

Looking ahead, the outlook for this research is promising. The complexity of calculating the persistence transformation is quadratic to the number of critical points of the input. Furthermore, by treating low persistence features as noise, we can adapt the algorithm to filter out such peaks with a single hyper-parameter, leading to a denoising functionality of the persistence transformation and a more concise output.

In addition to the outlook presented earlier, exploring the characteristics of critical points offers promising avenues for further improving the algorithms and matchings. Applying Morse theory (\cite{milnor1963morse}) to analyze critical points could provide valuable insights and refinements to our approach. Additionally, leveraging Prof. Kozlov's work on optimizing matchings (see \cite{kozlov2020combinatorial}) could lead to significant enhancements in the overall concept.

By delving deeper into the nature of critical points and their interplay with the proposed methods, we have the potential to unlock novel techniques and strategies to advance topological data analysis. Such investigations will undoubtedly contribute to the robustness and versatility of our approach and pave the way for even more impactful applications in various domains.

Lastly, the adaptability of this work to specific applications allows for customized implementations that precisely address the unique demands of various domains. As we continue to refine and extend this method, we are confident in its ability to contribute significantly to advancing topological data analysis in diverse real-world scenarios.

\begin{backmatter}
\section*{Acknowledgements}
The authors extend their sincere gratitude to Prof. Dmitry Feichtner-Kozlov at the University of Bremen for supervising GK and LR. Also, we sincerely thank Lukas Mentz for his constructive criticism of the manuscript. The financial support by the German Research Foundation via the RTG 2224, titled "$\pi^{3}$: Parameter Identification - Analysis, Algorithms, Implementations" is gratefully acknowledged.

\section*{Authors' information} 
\begin{itemize}
    \item Gideon Klaila, klailag@uni-bremen.de, ORCID: 0009-0002-2861-2095
    \item Anastasios Stefanou, stefanou@uni-bremen.de, ORCID: 0000-0002-5408-9317
    \item Lena Ranke, lranke@uni-bremen.de, ORCID: 0009-0003-4258-1608
\end{itemize}

\section*{Authors' contributions}
This paper is based on joined work of the three main authors as equal contributors. Gideon Klaila developed the theoretical formalism, contributed the main conceptual ideas and performed the analytic calculations. Both Anastasios Stefanou and Lena Ranke verified the results and contributed to the final version of the manuscript while Anastasios Stefanou helped supervise the project.

\section*{Competing interests}
The authors declare that they have no competing interests.

\section*{Declarations}

\subsection*{Ethical Approval}
Not applicable.

\subsection*{Funding}
Deutsche Forschungsgemeinschaft. Grant Number: RTG 2224

\subsection*{Availability of data and materials }
Not applicable.

\bibliography{sn-bibliography}

\begin{thebibliography}{27}
\providecommand{\natexlab}[1]{#1}
\providecommand{\url}[1]{{#1}}
\providecommand{\urlprefix}{URL }
\providecommand{\doi}[1]{\url{https://doi.org/#1}}
\providecommand{\eprint}[2][]{\url{#2}}
 \bibcommenthead

\bibitem[{Aichler and Walch(2015)}]{aichler2015maldi}
Aichler M, Walch A (2015) Maldi imaging mass spectrometry: current frontiers
  and perspectives in pathology research and practice. Laboratory investigation
  95(4):422--431

\bibitem[{Atkeson and McIntyre(1986)}]{atkeson1986robot}
Atkeson C, McIntyre J (1986) Robot trajectory learning through practice. In:
  Proceedings. 1986 IEEE International Conference on Robotics and Automation,
  IEEE, pp 1737--1742

\bibitem[{Berwald et~al(2018)Berwald, Gottlieb, and
  Munch}]{berwald2018computing}
Berwald JJ, Gottlieb JM, Munch E (2018) Computing wasserstein distance for
  persistence diagrams on a quantum computer. arXiv preprint arXiv:180906433

\bibitem[{Bubenik et~al(2015)Bubenik, De~Silva, and Scott}]{bubenik2015metrics}
Bubenik P, De~Silva V, Scott J (2015) Metrics for generalized persistence
  modules. Foundations of Computational Mathematics 15:1501--1531

\bibitem[{Carlsson(2009)}]{carlsson2009topology}
Carlsson G (2009) Topology and data. Bulletin of the American Mathematical
  Society 46(2):255--308

\bibitem[{Chazal et~al(2009)Chazal, Cohen{-}Steiner, Glisse, Guibas, and
  Oudot}]{chazal2009stability}
Chazal F, Cohen{-}Steiner D, Glisse M, et~al (2009) Proximity of persistence
  modules and their diagrams. In: Hershberger J, Fogel E (eds) Proceedings of
  the 25th {ACM} Symposium on Computational Geometry, Aarhus, Denmark, June
  8-10, 2009. {ACM}, pp 237--246, \doi{10.1145/1542362.1542407},
  \urlprefix\url{https://doi.org/10.1145/1542362.1542407}

\bibitem[{Cohen-Steiner et~al(2005)Cohen-Steiner, Edelsbrunner, and
  Harer}]{cohen2005stability}
Cohen-Steiner D, Edelsbrunner H, Harer J (2005) Stability of persistence
  diagrams. In: Proceedings of the twenty-first annual symposium on
  Computational geometry, pp 263--271

\bibitem[{Das et~al(2023)Das, Anand, and Chung}]{das2023topological}
Das S, Anand DV, Chung MK (2023) Topological data analysis of human brain
  networks through order statistics. Plos one 18(3):e0276419

\bibitem[{D{\l}otko and Gurnari(2022)}]{dlotko2022euler}
D{\l}otko P, Gurnari D (2022) Euler characteristic curves and profiles: a
  stable shape invariant for big data problems. arXiv preprint arXiv:221201666

\bibitem[{Edelsbrunner et~al(2002)Edelsbrunner, Letscher, and
  Zomorodian}]{edelsbrunner2002topological}
Edelsbrunner, Letscher, Zomorodian (2002) Topological persistence and
  simplification. Discrete \& Computational Geometry 28:511--533

\bibitem[{Edelsbrunner and Harer(2010)}]{Edelsbrunner2010}
Edelsbrunner H, Harer JL (2010) Computational topology: an introduction.
  American Mathematical Society, Providence, USA

\bibitem[{Edelsbrunner et~al(2008)Edelsbrunner, Harer
  et~al}]{edelsbrunner2008persistent}
Edelsbrunner H, Harer J, et~al (2008) Persistent homology-a survey.
  Contemporary mathematics 453(26):257--282

\bibitem[{Ge et~al(2010)Ge, Zheng, Hao, Shao, Wang, and
  Luterbacher}]{ge2010temperature}
Ge QS, Zheng JY, Hao ZX, et~al (2010) Temperature variation through 2000 years
  in china: An uncertainty analysis of reconstruction and regional difference.
  Geophysical Research Letters 37(3)

\bibitem[{Ghrist(2008)}]{ghrist2008barcodes}
Ghrist R (2008) Barcodes: the persistent topology of data. Bulletin of the
  American Mathematical Society 45(1):61--75

\bibitem[{Klaila et~al(2023)Klaila, Vutov, and Stefanou}]{klaila2023supervised}
Klaila G, Vutov V, Stefanou A (2023) Supervised topological data analysis for
  maldi mass spectrometry imaging applications. BMC bioinformatics 24(1):279

\bibitem[{Koseki et~al(2023)Koseki, Hayashi, Kojima, Hirose, and
  Shimamura}]{koseki2023topological}
Koseki J, Hayashi S, Kojima Y, et~al (2023) Topological data analysis of
  protein structure and inter/intra-molecular interaction changes attributable
  to amino acid mutations. Computational and Structural Biotechnology Journal
  21:2950--2959

\bibitem[{Kozlov(2020)}]{kozlov2020combinatorial}
Kozlov DN (2020) A combinatorial method to compute explicit homology cycles
  using discrete morse theory. Journal of Applied and Computational Topology
  4(1):79--100

\bibitem[{Kudryavtseva(2009)}]{kudryavtseva2009uniform}
Kudryavtseva EA (2009) Uniform morse lemma and isotopy criterion for morse
  functions on surfaces. Moscow University Mathematics Bulletin 64(4):150--158

\bibitem[{M{\'e}moli(2011)}]{memoli2011gromov}
M{\'e}moli F (2011) Gromov--wasserstein distances and the metric approach to
  object matching. Foundations of computational mathematics 11:417--487

\bibitem[{Milnor(1963)}]{milnor1963morse}
Milnor JW (1963) Morse theory. 51, Princeton university press

\bibitem[{Nicolau et~al(2011)Nicolau, Levine, and
  Carlsson}]{nicolau2011topology}
Nicolau M, Levine AJ, Carlsson G (2011) Topology based data analysis identifies
  a subgroup of breast cancers with a unique mutational profile and excellent
  survival. Proceedings of the National Academy of Sciences 108(17):7265--7270

\bibitem[{Otter et~al(2017)Otter, Porter, Tillmann, Grindrod, and
  Harrington}]{otter2017roadmap}
Otter N, Porter MA, Tillmann U, et~al (2017) A roadmap for the computation of
  persistent homology. EPJ Data Science 6:1--38

\bibitem[{Radtke et~al(2016)Radtke, Schwab, Wolf, Freitag, Alt, Kesch,
  Popeneciu, Huettenbrink, Gasch, Klein et~al}]{radtke2016multiparametric}
Radtke JP, Schwab C, Wolf MB, et~al (2016) Multiparametric magnetic resonance
  imaging (mri) and mri--transrectal ultrasound fusion biopsy for index tumor
  detection: correlation with radical prostatectomy specimen. European urology
  70(5):846--853

\bibitem[{Silva et~al(2018)Silva, Munch, and
  Stefanou}]{stefanou2018interleaving}
Silva V, Munch E, Stefanou A (2018) Theory of interleavings on categories with
  a flow. Theory and Applications of Categories 33:583--607

\bibitem[{Singh et~al(2023)Singh, Farrelly, Hathaway, Leiner, Jagtap, Carlsson,
  and Erickson}]{singh2023topological}
Singh Y, Farrelly CM, Hathaway QA, et~al (2023) Topological data analysis in
  medical imaging: current state of the art. Insights into Imaging 14(1):1--10

\bibitem[{Ver~Hoef et~al(2023)Ver~Hoef, Adams, King, and
  Ebert-Uphoff}]{ver2023primer}
Ver~Hoef L, Adams H, King EJ, et~al (2023) A primer on topological data
  analysis to support image analysis tasks in environmental science. Artificial
  Intelligence for the Earth Systems 2(1):e220039

\bibitem[{Yu and You(2021)}]{yu2021shape}
Yu B, You K (2021) Shape-preserving dimensionality reduction: An algorithm and
  measures of topological equivalence. arXiv preprint arXiv:210602096

\end{thebibliography}

\end{backmatter}

\end{document}